\documentclass[12pt]{article}

\usepackage{amsmath}

\usepackage{amssymb}

\usepackage{amscd}

\usepackage{bbm}

\newtheorem{theorem}{Theorem}[section]
\newtheorem{lemma}[theorem]{Lemma}

\newtheorem{proposition}[theorem]{Proposition}

\newtheorem{exa}[theorem]{Example}

\newtheorem{exas}[theorem]{Examples}

\newtheorem{prope}[theorem]{Property}

\newtheorem{defini}[theorem]{Definition}
\newenvironment{definition}{\begin{defini} \em}{\end{defini}}

\newtheorem{rema}[theorem]{Remark}

\newenvironment{equationth}{\stepcounter{theorem}\begin{equation}}{\end{equation}}

\def\ind{{\rm{Ind}}}

\def\BC{{\Bbb C}}

\def\Z{\mathbbm{Z}}

\def\R{\mathbbm{R}}

\def\bS{\mathbb{S}}

\def\C{\mathbbm{C}}

\def\B{\mathbbm{B}}

\def\a{\alpha}

\def\b{\beta }

\def\e{\varepsilon}

\def\O{\mathcal O}

\def\v {\vskip.1cm}

\makeindex

\begin{document}

\title {\bf  \Large 
{INDICES OF  VECTOR FIELDS  ON  SINGULAR VARIETIES: AN OVERVIEW}}


\author{\Large  {Jos\'e Seade}\thanks{Supported by CNRS (France), CONACYT and DGAPA-UNAM (Mexico)} \\  \\
 }


\date {\it{$\qquad \qquad   \qquad \qquad\qquad \qquad   \qquad \qquad $ Dedicado a Jean-Paul, \\
$\qquad \qquad   \qquad \qquad\qquad \qquad   \qquad \qquad$ con gran respeto y  afecto} }

\setcounter{section}{-1}




\maketitle

\section{Introduction}

The Poincar\'e-Hopf total index of a vector field with isolated singularities on a smooth, closed manifold $M$  can be regarded as
the obstruction to constructing a non-zero section of the tangent bundle $TM$. In this way it
extends naturally to complex vector bundles in general and leads to the notion of Chern classes.
When working with singular analytic varieties, it is thus natural to ask what should be
the notion of ``the index" of a vector field. 

Indices of vector fields on singular varieties were first considered by M. H. Schwartz 
in \cite {Sch1, Sch2}
in her study of  Chern classes for singular varieties.  For her purpose
there was no point in considering vector fields in general, but only a special
class of vector fields (and frames) that she called ``radial'', which are obtained by the important
process of {\it  radial extension}.  
The generalisation of this index to other vector fields  was defined independently in 
\cite {KT, EG0, SS3} (see also \cite {ASV}), and its extension 
for frames in general was done in \cite {BLSS2}. This index, that we
call {\it  Schwartz index}, is 
sometimes called ``radial index'' because it measures how far  the vector field is from being radial. In \cite
{EG0, ASV} this index is defined also for vector fields on real analytic varieties.

MacPherson in \cite {MP} introduced  the local Euler obstruction, also for 
constructing Chern classes of singular complex algebraic varieties.
In \cite {BS} this invariant was defined via vector fields, interpretation
that was essential  to prove (also in \cite {BS}) that the
Schwartz classes of a singular variety coincide with  MacPherson's classes. This
 viewpoint brings the local Euler obstruction into the frame-work of ``indices of
vector fields on singular varieties'' and yields to another index, 
that we may call {\it  the local Euler obstruction }
of the vector field at each isolated singularity; the Euler obstruction of the singular variety 
corresponding to the case of the radial vector field. This index relates to the previously mentioned 
Schwartz index by a formula known as the ``Proportionality Theorem'' of \cite {BS}. When
the vector field is determined by the gradient of a function on the singular variety, this 
local Euler obstruction is the defect studied in \cite {BMPS}.

On the other hand, one of the basic properties of the local index of Poincar\'e-Hopf is that it is 
stable under perturbations.
In other words, if $v$ is a vector field on an open set $U$ in $\R^n$  and $x \in U$ is an
isolated singularity of $v$, and if we perturb $v$ slightly, then its singularity at $x$ may split into several
 singular points of the new vector field  $\hat v$, but the sum of the indices of $\hat v$ at
these singular points equals the index of $v$ at $x$. If we now consider an analytic variety $V$ 
defined, say, by a holomorphic function $f:  (\C^n,0) \to (\C,0)$ with an isolated critical
 point at $0$, and if $v$ is a vector
field on $V$, non-singular away from $0$, then one would like ``the index" of $v$ at $0$ to be stable
under small perturbations of both, the function $f$ and the vector field $v$. 
The extension of this index
to the case of vector fields on isolated complete intersection singularity germs (ICIS for short) 
is immediate.
This leads naturally to another
concept of index, now called the {\it GSV-index}, 
introduced in \cite{Se1, GSV, SS3}.  There is also the analogous index for continuous vector fields on real analytic varieties
(see \cite {ASV, GMr1, GMr2}). 

One also has the virtual
index, introduced in \cite {LSS} for holomorphic
vector fields; the extension to continuous vector fields is immediate and
was done in \cite {SS3, BLSS2}. This index  is defined via  Chern-Weil theory. The idea is that
the usual Poincar\'e-hopf index can be regarded as a localisation, at the singular points of a vector field,
of the $n^{th}$-Chern class of a manifold. Similarly, for an ICIS $(V,0)$ in $\C^{n+k}$, defined by functions
$f= (f_1,\cdots,f_k)$, one has a localisation at $0$ of the top Chern class of the 
ambient space, defined by the gradient vector fields of the $f_i$ and the given vector field, tangent to $V$.
This localisation defines {\it the  virtual index} of the vector field;  this definition extends to a rather 
general setting, providing a topological way for looking at the top Chern class of the so-called
virtual tangent bundle of singular varieties which are local complete intersections.
In the case envisaged above, when $(V,0)$ is an  ICIS, this index coincides with the 
GSV-index.

Another remarkable property of the local
 index of Poincar\'e-Hopf is that in the case of germs of 
 holomorphic vector fields in $\BC^n$ with an isolated singularity at $0$, the local 
 index equals the integer:
 \[\hbox {dim} \;\O_{\BC^n,0} \big/ (X_1,\cdots,X_n)\,, \tag{*}\]
 where $(X_1,\cdots,X_n)$ is the ideal generated by the components of the vector field. 
This and other facts motivated the search for algebraic formulae for the index of  vector fields
on singular varieties.  The {\it homological index} of Gomez-Mont \cite {Gom} 
is a beautiful answer to that search.
It considers an isolated singularity germ $(V,0)$ of any dimension, and a holomorphic vector field on 
$V$, singular only at $0$. One has the K\"ahler differentials on $V$,
 and a Koszul complex $(\Omega^\bullet_{V,0}, v)$:
\[
0 \to \Omega^n_{V,0}
\to \Omega^{n-1}_{V,0} \to ... \to \O_{V,0}\to 0\,,
\]
where the arrows are given by contracting forms by the vector field $v$.
The homological index of $v$ is defined to be the Euler characteristic of this complex. 
When the ambient 
space $V$ is smooth at $0$, the complex is exact in all dimensions, except in degree $0$ where the 
corresponding homology group  has dimension equal to the local
 index of Poincar\'e-Hopf of $v$. When $(V,0)$ is a hypersurface germ, this index coincides with the 
GSV-index,
but for more general singularities the homological index is still waiting to be understood!

In fact, in \cite {EGS} there is given the corresponding notion of 
{\it homological index} for holomorphic 1-forms on singular varieties, and recent work of Sch\"urmann
 throws light into this, yet mysterious, invariant.

When considering smooth (real) manifolds, the tangent and cotangent bundles are canonically
isomorphic
and it does not make much difference to consider either vector fields or 1-forms in order to
define their indices and their relations with  characteristic classes.
 When the ambient space is a complex manifold, this is no longer the case, but there are still
ways for comparing indices of vector fields and 1-forms, and to use these to study Chern classes
of manifolds. To some extent this is also  true  for singular varieties,  but there are however
important differences and each of the two settings has its own advantages.

The first time that indices of  1-forms on singular varieties appeared in the literature
 was in MacPherson's work \cite {MP}, where he defined the local Euler obstruction in this way.
But the systematic study of these indices was begun by  W. Ebeling and S. Gusein-Zade  in 
a series of articles (see for instance \cite{EG2, EG3, EG4, EG5}). This has been, to some extent, a 
study parallel to the one for vector fields, outlined in this article. Also along this lines is 
\cite{BSS2}, which adapts to 1-forms  the radial extension technique of M. H. Schwartz and proves the
corresponding Proportionality Theorem.

Also,  J. Sch\"urmann in his book \cite {Schu2} introduces  powerful methods to studying singular
varieties via micro-local analysis and Lagrangian cycles, and much of the theory of indices of 1-forms 
can also be seen in that way. Furthermore, he has recently found a remarkable method 
for assigning an index of 1-forms to each constructible function on a Whitney stratified 
complex analytic space,
 in such a way that each of the known indices corresponds to a particular choice of a 
constructible function. This is closely related to 
MacPherson work in \cite {MP} for defining characteristic classes of singular
varieties.

In this article
we  briefly review the various  indices of vector fields on singular varieties. I am presently
 working  with Jean-Paul Brasselet and Tatsuo Suwa writing  \cite {BSS3}, a monograph with a detailed account 
of all these indices, through the viewpoints of algebraic topology (obstruction theory) and 
differential geometry (Chern-Weil theory), together with their relations with Chern classes of singular varieties. This will include some applications of 
these indices to other fields of singularity theory.

This article grew from my talk in the singularities meeting at the CIRM in Luminy
in celebration of the 60th anniversary of Jean-Paul Brasselet, 
and I want to thank 
the organizers for the invitations to participate in that meeting and to write these notes, particularly to 
Anne Pichon. I am also grateful to Tatsuo Suwa, Jean-Paul Brasselet and
J\"org Sch\"urmann for many helpful conversations.

\section{The Schwartz index}

 Consider first 
 the case when the ambient space is an affine irreducible complex analytic variety 
$V \subset \BC^N$ of  dimension  $n > 1$ with an isolated singularity at $0$. Let $U$ be an open ball 
around 
$0 \in \BC^N$, small enough so that every sphere in $U$ centered at $0$ meets $V$ transversally 
(see \cite {Mi2}).
For simplicity we restrict the discussion to $U$ and  set $V = V \cap U$. 
Let $v_{rad}$ be a continuous vector field on $V \setminus \{0\}$ which is transversal 
(outwards-pointing) to all spheres around $0$, and scale it so that it extends to a continuous section of 
$T \BC^N|_{V}$ with an isolated zero at $0$.  We call $v_{rad}$ {\it a radial vector field} 
at $0 \in V$. 
 Notice  $v_{rad}$  can be  further extended to a {\it radial} 
 vector field $v_{rad}^\#$  on all of $U$, {\it i.e.}  
transversal to all spheres centered at $0$.
By definition {\it the Schwartz index} of
$v_{rad}$ is the Poincar\'e-Hopf index at $0$ of the radial extension $v_{rad}^\#  $, so it is
$1$.  Of course we could have started with the zero-vector at $0$, then 
extend this to $v_{rad}$ on $V$  as above, and then extend it further to all of $U$ being transversal to all 
the spheres, getting the same answer;  this is the viewpoint that generalises when the 
singular set of $V$ has dimension more than $0$. 

Let us continue with the case when $V$ has an isolated singularity at $0$, and assume now that $v$ is
a continuous vector field on $V$ with an isolated singularity at $0$. By this we mean a continuous 
section 
$v$ of $T\BC^N|_V$ which is tangent to $V^* = V \setminus \{0\}$.
We want to define {\it the Schwartz index} of $v$; this index somehow  
measures the ``radiality'' of the vector field. It has various names in the 
literature
(c.f. \cite {KT,SS3,EG1,ASV}), one of them being {\it radial index}.

Let $ v_{rad}$ be a radial 
vector field 
at $0$, {\it i.e.}  $ v_{rad}$ is transversal, outwards-pointing, to the intersection of $V$ 
with every 
sufficiently small sphere $\bS_{\e}$ centered at $0$. We may now define the difference
between $v$ and $v_{rad}$ at $0$: consider small spheres 
$\bS_{\e}$,
$\bS_{\e'}$; $ \e > \e' > 0$, and let $w$ be a vector field on the cylinder $X$ in 
$V$ bounded 
by the links $K_{\e} = \bS_{\e} \cap V$ and $K_{\e'} = \bS_{\e'} \cap V$, such 
that $w$
 has finitely many singularities in the interior of $X$, it restricts to $v$ on 
$K_{\e}$ and 
to $v_{rad}$ on $K_{\e'}$. The {\it difference} 
$d(v, v_{rad}) = d(v, v_{rad}; V)$ of $v$ and $ v_{rad}$ is:

$$d(v, v_{rad}) = \ind_{PH}(w;X)\,,$$
the Poincar\'e-Hopf index of $w$ on $X$. Then define {\it the Schwartz} (or radial)
{\it index} of $v$ at $0 \in V$  to be:
$$ \ind_{Sch}(v,0;V) = 1 + d(v, v_{rad})\,.$$

The following result is well known (see for instance \cite {KT, EG1, SS3, ASV}). For vector fields 
with radial singularities, this is a special case of the work of 
M. H. Schwartz; the general case follows easily from this.

\begin{theorem}\label{sch1} Let $V$ be a compact complex analytic variety with isolated singularities $q_1,\cdots,q_r$
 in a complex manifold $M$, and let $v$ be a continuous vector field on $V$, singular at the $q_i\,'s$ and 
possibly at some other isolated points in $V$. Let $\ind_{Sch}(v;V)$ be the sum of the Schwartz indices of $v$
at the $q_i$ plus its Poincar\'e-Hopf index at the singularities of $v$ in the regular part of $V$. Then:
$$\ind_{Sch}(v;V) = \chi(V)\;.$$
\end{theorem}

The proof is fairly simple and we refer to the literature for details.

\vskip.2cm
The idea for defining  the Schwartz index in general, when the singular set has dimension more 
than $0$, 
 is similar in spirit to the case above, but it presents some technical difficulties. 
Consider a compact, complex analytic variety $V$ of dimension $n$ embedded in a complex manifold $M$,
equipped with a Whitney stratification $\{V_\a\}_{\a \in A}$ adapted to $V$. The starting point 
to define the Schwartz index of a vector field is the {\it radial extension} introduced by M. H. 
Schwartz. To  explain this briefly, let $v$ be a vector field defined on a neighbourhood of $0$ in the 
 stratum $V_\a$  of $V$ that contains $0$. 
The fact that the stratification is Whitney implies 
(see \cite {Sch1, Sch2, Br2} 
for details) that one can make a {\it parallel extension} of $v$ to a stratified vector field
$v'$ on a neighbourhood of $0$ in $M$. Now, if $0$ is an isolated singularity of $v$ on $V_\a$, then
$v'$ will be singular in a disc of dimension 
$(\hbox{dim}_\R M - \hbox{dim}_\R\, V_\a)$, transversal to  $V_\a$ in 
$M$ at $0$. So this extension is not good enough by itself. We must add to it another vector field
$v''$: the gradient of the square of the function ``distance to'' $V_\a$, defined near $0$. This vector field is 
transversal to the boundaries of all tubular neighbourhoods of $V_\a$ in $M$; using the Whitney
conditions we can make $v''$ be a continuous, stratified vector field near $0$. The zeroes of $v''$ are the points in $V_\a$. Adding $v'$ and $v''$ at each point near $0$ we get a stratified, 
continuous vector field $v^\#$ defined on a neighbourhood of $0$ in $M$, which restricts to the given
vector field $v$ on $V_\a$. This vector field has the additional property of being radial in 
all directions which are normal to the stratum $V_\a$. In other words, if we
take a small smooth disc $\Sigma$ in $M$ transversal to $V_\a$ at $0$ of 
dimension complementary to that of $V_\a$. Then the restriction of $v^\#$ to
 $\Sigma$ can be projected into a vector field tangent to  $\Sigma$ with 
Poincar\'e-Hopf index $1$ at $0$. Hence the Poincar\'e-Hopf index of
$v$ on the stratum $V_\a$ equals the Poincar\'e-Hopf index of $v^\#$ 
in the ambient space $M$: this is a basic
property of the vector fields obtained by radial extension.

\begin{definition} The {\it Schwartz index} of $v$ at $0 \in V_\a \subset V$ is defined to be the
 Poincar\'e-Hopf index at $0$ of its radial extension $v^\#$ to a neighbourhood of $0$ in $M$. 
\end{definition}

From the previous discussion we deduce:
\begin{proposition}
If the stratum $V_\a$ has dimension $> 0$, the Schwartz index  of $v$ 
equals the  Poincar\'e-Hopf index of $v$  at $0$ regarded as a vector field 
on the stratum $V_\a$.
\end{proposition}

Now, more generally, let $v$ be a  stratified vector field  on $V$ with an isolated 
singularity at $0 \in V \subset M$. 
 Let  $ v_{rad}$ be a stratified radial vector field 
at $0$, {\it i.e.}  $v_{rad} $ is transversal (outwards-pointing) to the intersection of $V$ 
with every 
sufficiently small sphere $\bS_{\e}$ in $M$ centered at $0$, and it is tangent to each stratum.
We  define the difference between $v$ and $v_{rad}$ at $0$ as follows. 
Consider sufficiently small spheres $\bS_{\e}$,
$\bS_{\e'}$ in $M$, $ \e > \e' > 0$, and put the vector field $v$ on 
$K_{\e} = \bS_{\e} \cap V$ and  $v_{rad}$ on $K_{\e'} = \bS_{\e'} \cap V$. We now use the Schwartz's technique of radial extension explained 
before, to get  a stratified 
 vector field $w$  on the cylinder $X$ in $V$ bounded 
by the links $K_{\e}$ and $K_{\e'} $, 
such 
that $w$ extends $v$ and  $v_{rad}$, it has finitely many singularities in the interior of $X$ and
at each of these singular points its index in the stratum equals its index in the ambient space 
$M$ (see \cite {Br2} for details). 
The {\it difference} of $v$ and $v_{rad}$ is defined as:
$$d(v,v_{rad}) = \sum \ind_{PH}(w;X)\,,$$
where the sum on the right runs over the singular points of $w$ in $X$ and each singularity is
being counted with the local Poincar\'e-Hopf index of $w$ in the corresponding stratum. As in the work 
of M. H. Schwartz, we can check that this integer does not depend on the choice of $w$.

\begin{definition}\label{def-radial} {\it The Schwartz} (or radial)
{\it index} of $v$ at $0 \in V$  is:
$$ \ind_{Sch}(v,0;V) = 1 + d(v,v_{rad})\,.$$
\end{definition}

It is clear that if $V$ is smooth at $0$ then this index coincides with the usual
 Poincar\'e-Hopf index; it also coincides with the index defined above when $0$ is an 
isolated singularity of $V$ and with the usual  index of M. H. Schwartz 
  for vector fields obtained by radial 
extension. In order to give a unified
 picture of what this index measures in the various cases, it is useful to introduce a  concept
 that picks up one of the essential 
properties of the vector fields obtained by  radial extension:

\begin{definition}\label{normally radial}\index{normally radial}
A stratified vector field on $V$ 
 is {\it normally radial}  at  $0 \in V_\a$ if it is radial in the direction of each stratum
$V_\b \ne V_\a$ containing $0$ in its closure. 
\end{definition}

In other words, $v$ is normally radial if its projection to each small disc 
$\Sigma$ around $0$, which is transversal to $V_\a$ at $0$ and has dimension
$(\hbox{dim}_\R\,M \,-\, \hbox{dim}_\R\, V_\a)$,
 is a radial vector field in $\Sigma$, 
{\it i.e.} it is transversal to each sphere in $\Sigma$ centered at $0$. The 
vector fields obtained by  radial extension 
satisfy this condition at all points.

The proof of the following proposition is immediate from the definitions.

\begin{proposition}\label{Schwartz unified}
Let $v$ be a  stratified vector field on $V$ with an isolated singularity at  $0$,
and let  $ V_\a$ be the Whitney stratum that contains $0$. If  $v$ is normally radial at 
$0$, then its Schwartz index $ \ind_{Sch}(v,0;V)$  
equals its Poincar\'e-Hopf index $\ind_{PH}(v,0;V_\a)$ in  $ V_\a$. Otherwise,  its 
 Schwartz index $ \ind_{Sch}(v,0;V)$ is the sum:
$$ \ind_{Sch}(v,0;V) \,=\, \ind_{PH}(v,0;V_\a) \,+\, \sum_{\b \ne \a} d(v, v_{rad}; V_\b)\,,$$
where $\ind_{PH}(v,0;V_\a)$ is defined to be $1$ if the stratum $V_\a$ has dimension $0$, and the sum
in the right runs over all strata that contain $V_\a$ in their closures; $d(v, v_{rad}; V_\b)$ is the difference in each stratum  $ V_\b$
 between $v$ and a stratified radial vector field $v_{rad}$ at $0$ .
\end{proposition}

\section{The local Euler obstruction}

Let  $ (V,0) $  be a reduced, pure-dimensional complex analytic singularity 
germ of dimension $n$ in an open set $ U \subset \C^N$. Let $G(n,N)$ denote the
Grassmanian of complex $n$-planes in $\C^N$. On the regular part $V_{reg}$ of
$V$ there is a map $\sigma : V_{ reg} \to U\times G(n,N)$ defined by
$\sigma(x) = (x,T_x(V_{ reg}))$. The {\it Nash transformation} 
(or {\it Nash blow up}) $\widetilde   V$ of $V$ is
the closure of Im$(\sigma)$ in $ U\times G(n,N)$. It is a  (usually singular)
complex analytic space endowed with an analytic projection map 
$$\nu : \widetilde   V \to V$$ 
which is a biholomorphism away from $\nu^{-1}(Sing(V))\,,$  
where $Sing(V):=V-V_{reg}$. Notice each point  $y \in Sing(V) $ is being replaced
by all limits of planes $T_{x_i}V_{reg}$ for sequences $\{x_i\}$ in $V_{reg}$ converging to
$x$.

Let us denote by  $ U(n,N)$ the tautological bundle over
$G(n,N)$ and denote by ${\Bbb U}$ the corresponding trivial extension bundle
over $ U \times G(n,N)$. We denote by $\pi$ the projection map of this bundle.
Let $\widetilde   T$ be the restriction of  ${\Bbb U}$ to $\widetilde   V$, with
projection map $\pi$.  The bundle $\widetilde   T$ on $\widetilde   V$ is called 
{\it the  Nash bundle} of $V$.  
An element of $\widetilde   T$ is written $(x,T,v)$ where $x\in
U$, $T$ is a $d$-plane in $\C^N$ and $v$ is a vector in $T$. We have maps:
$$ \widetilde   T \buildrel{\pi} \over {\longrightarrow} \widetilde   V \buildrel{\nu} 
\over {\longrightarrow} V \,,
$$ 
where $\pi$ is the projection map of the Nash bundle over the Nash blow up 
$\widetilde   V$.

 Let us consider a complex analytic stratification $(V_\a)_{\a\in A}$ of $V$
satisfying the Whitney conditions. Adding the stratum $U\setminus V$ we
obtain a Whitney stratification of $U$.  Let us denote by $TU|_V$ the
restriction to $V$ of the tangent bundle of $U$. We know that a stratified vector field
$v$ on $V$ means a continuous section of $TU|_V$ such that if $x \in V_\a\cap
V$ then $v(x) \in T_x(V_\a)$. The Whitney condition (a) implies that given 
$x \in Sing(V)$, any limit $\mathcal T$ of tangent spaces of points in $V_{reg} = V - Sing(V)$
converging to $x$ contains the tangent space $T_xV_\a$ where $V_\a$ is 
the stratum that contains $x$. Hence one has the following lemma of \cite {BS}:

\begin{lemma}\label{lifting} Every stratified vector field $v$ on a set $A \subset V $
has a canonical lifting to a section $\widetilde   v$ of the Nash bundle $\widetilde   T$ over
 $\nu^{-1}(A) \subset \widetilde   V$.
\end{lemma}

Now consider a stratified  radial vector field $v(x)$ in a neighborhood of $\{
0\}$ in $V$; {\it i.e.}  there is $\varepsilon_0$ such that for
every $0<\varepsilon \leq \varepsilon_0$, $v(x)$ is pointing outwards the
ball $\Bbb
B_\varepsilon$ over the boundary  
$V \cap \Bbb S_\varepsilon$ with $\Bbb S_\varepsilon:=\partial \Bbb
B_\varepsilon$. 
Recall that, essentially by the Theorem of Bertini-Sard (see \cite {Mi2}), for $\varepsilon$
small enough the spheres ${\Bbb S}_\varepsilon$ are transverse to the strata
$(V_\a)_{\a \in A}$.

One has the following interpretation of the local Euler obstruction 
\cite {BS}. We refer to \cite {MP} for the original definition which uses
1-forms instead of vector fields.

\begin{definition}\index{Euler obstruction}\label{def. Euler obstruction}
Let $v$ be a stratified radial vector field on $V \cap \Bbb
S_\varepsilon$ and $\widetilde   v$ the lifting of $v$ on $\nu^{-1} (V\cap \Bbb
S_\varepsilon)$ to a section of the Nash bundle. The {\it local Euler obstruction} (or simply the Euler
obstruction) Eu$_V(0)$ is defined to be the obstruction to  extending  $\widetilde  
v$ as a nowhere zero section of $\widetilde   T$ over $\nu^{-1} (V\cap \Bbb
B_\varepsilon)$.
\end{definition}

More precisely, let ${\mathcal O} (\widetilde   v) \in H^{2d}\big(
\Bbb B_\varepsilon, \nu^{-1}(V \cap \Bbb S_\varepsilon)\big)\,$  
be
the
obstruction cocycle for  extending  $\widetilde   v$ as a nowhere zero section
of
$\widetilde   T$ inside $\nu^{-1} (V\cap \Bbb B_\varepsilon)$, where $\Bbb B_\varepsilon$
is a small ball around $0$ and $\Bbb S_\varepsilon$ is its boundary.
The local Euler obstruction $Eu_V(0)$ is the evaluation of ${\mathcal O}
(\widetilde   v)$  on the fundamental class of the pair $\big(\nu^{-1}(V \cap
\Bbb B_\varepsilon), \nu^{-1}(V \cap \Bbb S_\varepsilon)\big)$.  
The
Euler obstruction is an integer.

The following result summarises some basic properties of the Euler obstruction:

\begin{theorem}\label{properties Euler} The  Euler obstruction satisfies:
\begin{enumerate}
\item $Eu_V(0) = 1$ if $0$ is a regular point of $V$;

\item   $Eu_{V \times V'}(0 \times 0') = Eu_V(0) \cdot Eu_{V'}(0')$;

\item  If $V$ is locally reducible at $0$ and $V_i$ are its irreducible
components, then $Eu_{V}(0) = \sum Eu_{V_i}(0)$;

\item  $Eu_V(0)$ is a constructible function on $V$, in fact it is constant
on Whitney strata.
\end{enumerate}
\end{theorem}

These statements are all contained in \cite {MP}, except for the second part of (iv)
which is not explicitly stated there 
and we refer to \cite {BS, LT1} for a detailed proof.

More generally, 
for every point $x\in V$, we will denote by
$V_\a(x)$ the stratum containing $x$. Now suppose $v$ is a stratified vector field
on a small disc  $\B_x$ around  $x \in V$, and $v$ has an isolated singularity at $x$. By 
\ref{lifting} we have  that $v$ can be lifted to a section $\widetilde v$
of the Nash bundle $\widetilde T$ of  $V$ over $\nu^{-1}(\B_x  \cap V)$ and $\widetilde v$
is never-zero on $\nu^{-1}(\partial \B_x  \cap V)$. The obstruction for extending $\widetilde v$
without singularity to the interior of $\nu^{-1}(\B_x  \cap V)$ is a
 cohomology class in 
$H^{2n}(\nu^{-1}(\B_x  \cap V), \nu^{-1}(\partial \B_x  \cap V))$;
evaluating this class in the fundamental cycle $[\B_x, \partial \B_x]$ one gets
  {\it an index}  $Eu(v,x;V) \in \Z$ of $v$ at $x$.
If $v$ is radial at $x$ 
then $Eu(v,x;V)$ is by definition the local Euler obstruction of
$V$ at $x$, $Eu_V(x)$.

\begin{definition}
The integer  $Eu(v,x;V)$ is {\it the (local) Euler obstruction  of the stratified vector field} $v$ at $x \in V$.
\end{definition}


As mentioned in the introduction, this index is related to the Schwartz index by the Proportionality 
Theorem of \cite {BS}. To state this result,
 recall that  we introduced in section 1 the concept of
normally radial vector fields, which essentially characterises the 
 vector fields obtained by radial extension.

\begin{theorem}{\rm (Proportionality Theorem \cite {BS})}\label{proportionality-Eu}\index
{Proportionality theorem ! for Euler obstruction}
Let  $v$ be a stratified vector field on $V$ which is normally radial at 
a singularity $0 \in V_\a$. 
Then one has:
$$Eu(v,0;V) = \ind_{Sch}(v,0) \cdot  Eu_V(0)$$
where $ Eu_V(0)$ is the Euler obstruction of $V$ at $0$ and
$\ind_{Sch}(v,0)$ is the Schwartz index of $v$ at $0$.
\end{theorem}

In short, this theorem says that the obstruction $Eu(v,0;V)$ to extend the lifting 
$\widetilde   v$ as a section of the Nash bundle inside  $\nu^{-1}(V \cap \B_{\e}(p) $
is proportional to the Schwartz index of $v$ at $0$, the proportionality factor
 being precisely the local Euler obstruction. We refer to \cite {BSS4} for a short proof of this theorem.

\medskip
The invariant $Eu(v,0;V)$ was studied in \cite {BMPS} when $v$ is the ``gradient  vector field'' 
$\nabla_f$ of a function $f$ on $V$. More precisely, if $V$ has an isolated singularity at $0$
 and the (real or complex valued) differentiable function $f$ has an isolated critical point 
at $0$, then $v$ is truely the (complex conjugate if $f$ is complex valued) gradient of the restriction of
$f$ to $V \setminus \{0\}$. In general, if  $V$ has a non-isolated singularity at $0$ but $f$
has an isolated critical point at $0$ (in the stratified sense \cite {GMP, Le2}), then
$v$ is obtained essentially by projecting the gradient vector field of $f$ to the tangent space of 
the strata in $V$, and then using the Whitney conditions to put these together in a continuous, 
stratified vector field. One may define this invariant even if $f$ has non-isolated critical points,
using intersections of characteristic cycles (see \cite {BMPS}), and it is a measure of how far the germ $(V,0)$ is from
satisfying the local Euler condition (in bivariant theory) with respect to the function $f$. Thus 
it was called in \cite {BMPS} {\it the Euler defect} of $f$ at $(V,0)$. The Euler obstruction of
 MacPherson corresponds to the case when $f$ is the function distance to $0$.
As noticed in \cite {EG5}, this invariant can be also defined using the 1-form
$df$ instead of the gradient vector field. This avoids several
 technical difficulties and is closer to MacPherson's original definition of the
local  Euler obstruction.

In \cite {STV} it is proved that if $f$ has an isolated critical point at $0 \in V$ 
(in the stratified sense), then its ``defect'' equals the number of critical points in the regular 
part of $V$ of a morsification of $f$. This fact can also be deduced easily from \cite {Schu1}.

\section{The GSV-index}

Let us denote by $(V,0)$ the germ of a complex analytic $n$-dimensional, isolated
 complete intersection singularity, defined by a function 
$$ f = (f_1,...,f_k) \, \colon \, ({\C}^{n+k},0) \to \, ({\C}^k,0) \, ,$$
and let $v$ be a continuous vector field on $V$ singular only at
$0$. If $n=1$ we further assume (for the moment) that $V$ is irreducible.  
We use the notation of \cite {Lo1}: an ICIS means an isolated complete intersection singularity.

Since $0$ is an isolated singularity of $V$, it
follows that the (complex conjugate) gradient vector fields 
$\, \{\overline{\nabla} f_1,...,\overline{\nabla}
f_k\} \,$ are linearly independent everywhere away from $0$ and they are
normal to $V$.  Hence the set $ \, \{v, 
\overline{\nabla} f_1,...,\overline{\nabla}
f_k \} \,$ is a $(k+1)$-frame on $V\setminus \{0\}$.   Let $\,K = V \cap
\bS_{\varepsilon} \, $ be the link of $0$ in $V$. It is an oriented, 
real manifold of dimension $(2n-1)$ and the above  frame defines a map 
$$ \phi_v \, = \, (v,\overline{\nabla} f_1,...,\overline{\nabla} f_k ) \, 
\colon \, K \, 
\to W_{k+1}(n+k) \, ,$$
into the Stiefel manifold of complex $(k+1)$-frames in ${\C}^{n+k}$. Since 
$W_{k+1}(n+k)$ is 
simply connected, its first non-zero homology group is in dimension $(2n-1)$ and it
is isomorphic to $\Z$. Hence the map $ \phi_v$ has a well defined degree $\deg(\phi_v) \in \Z$. 
To define it we notice that $W_{k+1}(n+k)$ is a fibre bundle over $W_{k}(n+k)$ with fibre the sphere 
$\bS^{2n-1}$; if $(e_1,\cdots,e_{n+k})$ is the canonical basis of $\C^{n+k}$, then the fiber $\gamma$ over the
$k$-frame  $(e_1,\cdots,e_{k})$ determines the canonical generator $[\gamma]$ of $H_{(2n-1)}(W_{k}(n+k)) \cong \Z$.
If $[K]$ is the fundamental class of $K$, then  $ (\phi_v)_*[K] = \lambda \cdot [\gamma] $ for some integer
$\lambda$. Then the degree of $\phi_v$ is defined by:
$$ \deg(\phi_v) = \lambda\,.$$
Alternatively one can prove that every map from a 
closed oriented $(2n-1)$-manifold into $ \, W_{k+1}(n+k) \,$ factors by a map into
the fibre  $\gamma \cong  \bS^{2n-1}$,
essentially by transversality.  Hence $\phi_v$ represents an element in 
$\pi_{2n-1} W_{k+1}(n+k) \, \cong \, {\Z} \, ,$ so $\phi_v$ is classified by its 
degree.

\begin{definition}\label{GSV as degree}  The GSV-index of $v$ at $0 \in V$, $\ind_{GSV}(v,0;V)$,  is 
the degree 
of the above map $\phi_v$.
\end{definition}

This index depends not only on the topology of $V$ near $0$, but also on the way 
$V$ is embedded in the ambient space. For instance the singularities in 
${\C}^3$ defined by 
$$ \{x^2+y^7+z^{14} = 0\} \, \, \, \hbox {and} \, \, \{x^3+y^4+z^{12} = 0\} \, , $$ 
are orientation preserving homeomorphic, but one can prove that the GSV-index of the radial vector 
field is 79 in the first case and 67 in the latter; this follows from the fact (see \ref{GSV} below) that 
for radial vector fields the GSV-index is 
$1 + (-1)^{\hbox{dim}\, V} \mu$, where $\mu$ is the Milnor number, which in the examples 
above is known to be 78 and 66 respectively.

We recall that one has a Milnor fibration associated to the function
$f$, see \cite {Mi2, Ha, Lo1} and the Milnor fibre $F$ can be regarded as a compact 
$2n$-manifold with boundary $\partial F = K$.  Moreover, by the Transversal Isotopy 
Lemma (see for instance \cite {AR}) there is 
an ambient isotopy of the sphere $\bS_{\varepsilon}$ taking $K$ into $\partial F$, 
which can be extended to a collar of $K$, which goes into a collar of $\partial F$ 
in $F$.  Hence $v$ can be regarded as a non-singular vector field  on $\partial F$.

\begin{theorem}\label{GSV} This index has the following properties:

\noindent  
{ (i)}  The GSV-index of $v$ at $0$ equals the Poincar\'e-Hopf index of $v$ in 
the Milnor fibre:  $$\ind_{GSV}(v,0;V) = \ind_{PH}(v,F) \, .$$

\noindent
{ (ii)}  If $v$ is everywhere transversal to $K$, then
$$\ind_{GSV}(v,0;V) = 1 +(-1)^n \mu \, .$$

\noindent
{ (iii)}   One has: $$\ind_{GSV}(v,0;V) = \ind_{Sch}(v,0;V) + (-1)^n \mu \, ,$$
where $\mu$ is the Milnor number of $0$ and $\ind_{Sch}$ is the Schwartz index.
\end{theorem}

Notice that the last statement says that the Milnor number of $f$ equals (up to sign) the difference of the Schwartz and GSV indices of every vector field 
on $V$ with an isolated singularity (cf. \cite {EGS}).

In \cite {BSS1} there is a generalisation of this index to the case when the variety $V$ has 
non-isolated singularities, but the vector field is stratified and it has an isolated singularity.
In \cite {ASV} is studied the real analytic setting and relations with other invariants of real 
singularities are given.

If $V$ has dimension 1 and is not irreducible, then the GSV-index of vector fields on $V$ was 
actually introduced by M. Brunella in \cite {Bru1, Bru2} and by Khanedani-Suwa \cite {KS},
 in relation with the geometry of holomorphic
1-dimensional  foliations on complex surfaces. In this case one has two possible definitions of the 
index:  as the Poincar\'e-Hopf index of an extension of the vector field to a Milnor fibre, or 
as the sum of the degrees in \ref{GSV as degree} corresponding to the various branches of
$V$. One can prove \cite {Su, ASV} that for plane curves these integers differ by the intersection numbers of 
the branches of $V$.

\section{The Virtual Index}

We now let $V$ be a compact local complete intersection of
dimension $n$ in a manifold $M$ of dimension $m = n+k$,  defined
as the zero set of a holomorphic section $s$  of a holomorphic vector
bundle $E$ of rank $k$ over $M$. The singular set of $V$, $Sing(V)$,
 may have dimension $\ge 0$.
Let  $v$ be a $C^{\infty}$ vector field on $V$. We denote by
 $\Sigma$ the singular set of $v$, which is assumed to consist of
$Sing(V)$ and possibly some other connected components
in the regular part of $V$,
 disjoint from $Sing(V)$.

The virtual index is an  invariant that assigns an integer to each connected component $S$ of
$\Sigma$. When $S$ consists of one point, this index coincides with the GSV index, and for a 
component $S \subset V_{reg}$ this is just the sum of the local 
indices of the singularities into which $S$ splits under a morsification of $v$.

Given  a connected  component  $S$  of $\Sigma$,
  the idea  to define the virtual index $\ind_{Vir}(v;S)$ of $v$ at $S$ is to  
localize at $S$ a certain characteristic class. We know that if $V$ is smooth, then  the
usual Poincar\'e-Hopf local index can be regarded as the top Chern class $c_n(TD,v)$
of the tangent bundle of a small disc 
 around the singular point of the vector field, relative to
the section $v$ given by the vector field on the boundary of $D$. In other words, the
Poincar\'e-Hopf local index is obtained by localizing at  $0$ the top Chern class
$c_n(TD)$ using the vector field $v$. We can of course replace  the point $0$ by a
component  $S$ of  $\Sigma$ contained in $V_{reg}$; in this case we
 replace  $D$ by a compact tubular
neighbourhood $\mathcal  T$ of $S$ in $V$ and we localize $c_n(TV)|_{\mathcal  T}$ at $S$
using the vector field $v$, which is assumed to be
non-singular on $\mathcal  T \setminus S$. This means that we consider the Chern class
$c_n(TV|_{\mathcal  T})$ relative to $v$ on $\mathcal  T \setminus S$.
The class we get lives in
$H^{2n}(\mathcal  T, \mathcal  T \setminus S) \cong H_0(S) \cong \Z$. The integer that we get 
in this way
is the  Poincar\'e-Hopf index of $v$ at $S$, $\ind_{PH}(v;S)$, {\it i.e.}  the number of singularities
of a generic perturbation of $v$ inside   $\mathcal  T$, counted with signs.

The question now is what to do when $S$ is contained in the singular set of $V$, so 
there is not  a tangent bundle. The idea to define the virtual index is to make a
similar ``localisation'' using the vector field and the {\it virtual tangent bundle of $V$}, defined below.

To define this bundle we notice that
the restriction $E\vert_{V_{reg}}$
coincides with the (holomorphic) normal bundle  $N({V_{reg}})$ of the regular part $V_{reg} = V - Sing(V)$.
We denote by $TM$ the holomorphic  tangent bundle of $M$ and we set $N = E|_V$.

\begin{definition} (c.f. \cite {FJ})  The {\it virtual tangent bundle} of $V$ is
$$ \tau(V) = TM\vert_V-N\,,$$ regarded as an element in the complex K-theory $KU(V)$.
\end{definition}
\v

It is known that the equivalence class of this virtual bundle does not depend on the
choice of the embedding of $V$ in $M$.

We denote by
$$c_*(TM|_V) = 1 + c_1(TM|_V) .... + c_m(TM|_V)\,,$$
and
$$c_*(N) = 1 + c_1(N) .... + c_k(N)\,,$$
the total Chern classes of these bundles. These are elements in the cohomology ring
of $V$ and can be inverted, {\it i.e.}  there is a unique class $c_*(N)^{-1} \in H^*(V)$
such that 
$$c_*(N) \cdot c_*(N)^{-1} \,= \,c_*(N)^{-1} \cdot c_*(N) \,= \,1\,.$$ Using this one
has the total Chern class of the virtual tangent bundle defined in the usual way:
$$c_*(\tau(V)) = c_*(TM|_V) \cdot c_*(N)^{-1} \in H^*(V)\,.$$
The $i^{th}$ Chern class of $TM\vert_V-N$ is by definition the component of
$c_*(\tau(V) )$ in dimension $2i$, for $i = 1,...,n$.

It is clear that if $V$ is smooth, then its virtual tangent bundle is equivalent in
$KU(V)$ to its usual tangent bundle, and  the Chern classes of the virtual
tangent bundle are the usual Chern classes.

Consider the component $c_n(\tau(V) )$
 of $c_*(\tau(V))$ in dimension $2n$. This is the top Chern class of the virtual
tangent bundle. As we said before, the idea to define the local index of the vector
field $v$ at a component $S$ of $Sing(V)$ is to localize $c_n(\tau(V))$ at $S$ using $v$.
For this one needs to explain how to localize the Chern classes of the
virtual tangent bundle. This is carefully done in \cite {Su}, and we refer to that text 
 for a detailed account on the subject, particularly in relation with
indices of vector fields.

In the particular case when the component $S$ has dimension $0$, so that we can assume we
have a local ICIS germ $(V,0)$ of dimension $n$ in $\C^{n+k}$, defined by functions
$$f = (f_1,\cdots,f_k) : U \subset \C^{n+k} \to \C^k \,,$$ with $U$ an open set in $\C^{n+k}$, 
one has that 
the virtual tangent bundle of $V$ is:
$$\tau (V) \,=\, T\C^{n+k}|_V - (V \times T\C^k)\,.$$
If $\B$ denotes a small ball in $U$ around $0$, then one has the Chern class $c_{n+k}(T\B|_V)$ 
relative
to  the $(k+1)$-frame $(v, \overline \nabla f_1, \cdots, \overline \nabla f_k)$ on 
$\partial \B \cap V$. This is a cohomology class in 
$H^{2n+2k}(\B \cap V, \partial \B \cap V)$, and one can prove that 
its image in $\cong H_0(\B) \cong \Z$
under  the Alexander homomorphism is the virtual index of $v$ at $0$ (see \cite {LSS, SS3, Su}); which in this case coincides with the GSV-index.


\section{The Homological Index}
The basic references for  this section are the articles by Gomez-Mont and various 
co-authors, see \cite{Gom}  and also \cite{BGo, GG1, GG2,  GGM1, GGM2}. 
There are also  important algebraic formulas for the index of holomorphic 
vector fields (and 1-forms) given by various 
authors, as for  instance in \cite {LSS} (see also \cite {GSV, Kl1, Kl2, Kl3}). 
In the real analytic case, interesting  algebraic formulas for the index 
are given in \cite {EG0, GMr1, GMr2},  which generalize to singular hypersurfaces
the remarkable formula of Eisenbud-Levine and 
Khimshiashvili \cite {EL, Khi}, that expresses 
 the index of an analytic vector field in $\R^m$ as the signature of an appropriate bilinear form.
 Here we only describe (briefly) the homological index of 
holomorphic vector fields.

Let $(V,0)\subset (\BC^N, 0)$ be the germ of a complex analytic (reduced) variety of pure
dimension $n$ with   an isolated singular point at the origin.  A vector field $v$
  on $(V,0)$ can always be defined   as the restriction to $V$ of a vector field $\widetilde v$ in the 
ambient space which is tangent to $V \setminus \{0\}$; 
   $v$ is holomorphic if $\widetilde v$ can be chosen to be holomorphic. So we may write
  $v$ as $v = (v_1,\cdots, v_N)$ where the $v_i$ are restriction to $V$ of holomorphic functions on a 
neighbourhood of $0$ in $(\C^N, 0)$.
  
  It is worth noting that given every space $V$ as above, there are always holomorphic vector fields on
  $V$ with an isolated singularity at $0$. This (non-trivial) fact is indeed a weak 
form of stating  a stronger result (\cite[2.1, p. 19]{BGo}): in the space $\Theta(V,0)$ 
  of germs of holomorphic vector fields 
  on $V$ at $0$, those having an isolated singularity form a connected, dense open subset
  $\Theta_0(V,0)$.
   Essentially the same result implies also that every $v \in \Theta_0(V,0)$ can be extended
  to a germ of holomorphic vector field in $\C^N$ with an isolated singularity, though it can 
possibly be also 
 extended with a singular locus of dimension more that $0$, a fact that may be useful for explicit 
computations (c.f. \cite {Gom}).

 A (germ of) holomorphic $j$-form on $V$ at $0$ 
 means the restriction to $V$ of a holomorphic $j$-form on a neighbourhood of
 $0$ in $\C^N$; two such forms in $\C^N$ are equivalent if their restrictions to $V$
  coincide on a neighbourhood of $0 \in V$. 
  We denote by $\Omega^j_{V,0} $ the space of all such forms (germs); these are the K\"ahler 
differential
 forms on $V$ at $0$. 
  So, $\Omega^0_{V,0} $ is the local structure ring $\O_{(V,0)}$ 
  of holomorphic functions on $V$ at $0$, and each
  $\Omega^j_{V,0} $ is a $\Omega^0_{V,0} $ module. Notice that if the germ of $V$ at $0$ is 
  determined by $(f_1,\cdots,f_k)$ then one has:
  \[\Omega^j_{V,0} := \, \frac{\Omega^j_{\C^N,0} }{\big(f_1\Omega^j_{\C^N,0} + df_1 \wedge 
  \Omega^{j-1}_{\C^N,0}\,, \cdots, \, f_k\Omega^j_{\C^N,0} + df_k \wedge 
  \Omega^{j-1}_{\C^N,0}\big)} \;,
  \]
  where $d$ is the exterior derivative.

Now, given a holomorphic vector field $\widetilde v$ at $0 \in \C^N$ with an isolated singularity
at the origin, and a K\"ahler
 form $\omega \in \Omega^j_{\C^N,0} $, we can always contract $\omega$ by $v$ in the usual way, 
 thus getting a K\"ahler  form $i_v(\omega) \in \Omega^{j-1}_{\C^N,0} $. If $v = \widetilde v\vert_V$ 
is tangent to $V$, then contraction is well defined at the level of Kh\"aler
 forms on $V$ at $0$ and one gets a  
complex $(\Omega^\bullet_{V,0}, v)$:
\[
0 \to \Omega^n_{V,0} 
\to \Omega^{n-1}_{V,0} \to ... \to \O_{V,0}\to 0\,,
\]
where the arrows are contraction by $v$ and $n$ is the dimension of $V$; of course one also has Kh\"aler
 forms of degree $>n$, but those forms do not play a significant role here.
 We consider the homology groups of this complex:
$$
H_j(\Omega^\bullet_{V,0}, v) = \frac{Êker\,(\Omega^{j}_{V,0} \to 
\Omega^{j-1}_{V,0})}{Im\,(\Omega^{j+1}_{V,0} \to \Omega^j_{V,0}) }
$$
An important observation in \cite {Gom} is that if $V$ is regular at $0$, so that its germ at $0$ is
 that of $\C^n$ at the origin, 
 and if $v = (v_1,\cdots,v_n)$ has an isolated singularity 
at $0$, then this is the usual Koszul complex
(see for instance \cite[p.  688] {GH}), so that all its homology groups vanish for $j > 0$, while 
$$H_0(\Omega^\bullet_{V,0}, v)  \cong \O_{\C^n,0)}\big / (v_1,\cdots,v_n)\,.$$ 
In particular the complex is exact when $v(0) \ne 0$. Since 
the contraction maps are $\O_{V,0}$-modules maps, 
this implies that if $V$ has an isolated singularity at the origin, then the homology groups
of this complex are concentrated at $0$, and they are finite dimensional because the sheaves of Kh\"aler 
forms on $V$ are coherent. Hence, for $V$ a complex analytic affine germ with
 an isolated singularity at $0$ and
 $v$ a holomorphic vector field on $V$ with an isolated singularity at $0$,  
it makes sense to define:

\begin{definition}
  The {\it homological index}
$\,{\rm Ind}_{\rm hom}(v,0;V)$ of the holomorphic vector field $v$
on $(V, 0)$ is  the Euler characteristic of the above
complex:
\[
{\rm Ind}_{\rm hom}\,(v,0;V) = \sum_{i=0}^n
(-1)^{i}  h_i(\Omega^\bullet_{V,0},v)\,,
\]
where $h_i(\Omega^\bullet_{V,0},v)$ is the dimension of the
corresponding homology group as a vector space over $\C$.
\end{definition}

We recall that an important property of the Poincar\'e-Hopf local index is its stability under perturbations. 
This means that if we perturb $v$ slightly in a neighbourhood of an isolated singularity, then this zero of 
$v$ may split into a number of
isolated singularities of the new vector field $v'$, whose total number (counted with their local indices) 
is the index of $v$. When the ambient space $V$ has an isolated singularity at $0$, then every vector field
on $V$ necessarily vanishes at $0$, since in the ambient space the vector field defines a local 1-parameter family of diffeomorphisms.
 Hence every perturbation of $v$ producing a vector field tangent to $V$ must also vanish at $0$, but new 
singularities may arise with this perturbation.
The homological index also satisfies the stability under such perturbations. This is 
called the ``Law of Conservation of Number" in \cite {Gom, GG2}:

\begin{theorem}
 {\rm (Gomez-Mont \cite [Theorem 1.2]{Gom})} For every holomorphic vector field $v'$ 
on $V$ sufficiently close to $v$ one has: 
$${\rm Ind}_{\rm hom}\, (v,0;V)\,= \, {\rm Ind}_{\rm hom}\, (v',0;V) \, + 
\sum {\rm Ind}_{PH}(v')\;,$$ 
where ${\rm Ind}_{PH}$ is the local Poincar\'e-Hopf index and 
the sum on the right runs over the singularities of $v'$ at regular points of
$V$ near $0$. 
\end{theorem}

This result is a special case of a more general theorem in \cite {GG2}. 

\medskip
This theorem is a key property of the homological index. In particular this allows us to identify 
this index with the GSV-index when $(V,0)$ is a hypersurface germ \cite {Gom}. In fact, it is easy to see that 
the GSV-index also satisfies the above ``Law of Conservation of Number'' for vector fields on
 complete intersection germs. This implies that if both 
indices coincide for a given vector field on $(V,0)$, then they coincide for every 
 vector field on $(V,0)$, since the space $\Theta_0(V,0)$ is connected. 
Hence, in order to prove that both indices coincide for all vector fields 
on hypersurface (or complete intersection) germs, it is enough to show that given every such germ, 
there exists a holomorphic vector field $v$ for which the GSV and homological indices coincide. 
This is what Gomez-Mont does in \cite {Gom}. For that, he first gives a very nice algebraic formula
to compute the homological index of vector fields on hypersurface singularities, which he then uses to perform explicit computations and prove that, for
holomorphic vector fields on hypersurface singularities,
the homological index coincides with the GSV index. It is not known 
whether or not these indices coincide on complete intersection germs in general (c.f. \cite{EGS}).

\section{Relations with Chern classes of singular varieties}

The local index of Poincar\'e-Hopf  is the most basic invariant of a vector field at an isolated
singularity, and the theorem of Poincar\'e-Hopf about the total index of a vector field on a
manifold is a fundamental result, giving rise, in particular,  to obstruction theory and the theory
 of characteristic classes, such as the Chern classes of complex manifolds.

In the case of singular varieties, there are several definitions of characteristic classes,
given by various authors. Somehow they correspond to  the various extensions one has of the 
concept of ``tangent bundle'' as we go from manifolds to singular varieties, and they are closely related
to the indices of vector fields discussed above.

The first one is due to M.H. Schwartz in \cite {Sch1, Sch2}, considering a
 singular complex analytic variety $V$ embedded in a smooth complex manifold $M$ which is equipped with a
Whitney stratification adapted to $V$; she then replaces the tangent bundle by the union of 
tangent bundles of all the strata in $V$, and considers a class of stratified frames to define 
characteristic classes of $V$, which do not depend on $M$ nor on the various choices. 
These classes live in the cohomology 
groups of $M$ with support in $V$, {\it i.e.}  $H^*(M, M\setminus V; {\Bbb Z})$,
and they are equivalent to the usual Chern classes when $V$ is non-singular. The top degree Schwartz 
class is defined precisely using the Schwartz index presented in Section 1: consider a Whitney 
stratification of $M$ adapted to $V$, a triangulation $(K)$ compatible with the stratification, and
the dual cell decomposition $(D)$ (c.f. \cite {Br1, Br2} for details). By construction, the cells of 
$(D)$ are transverse to the strata. Consider now a 
stratified vector field $v$ on $V$ obtained by radial extension. Then (see \cite {Sch1, Sch2, BS, Br2})
the radial extension technique allows us to construct a vector field on a regular neighborhood
$U$ of $V$ in $M$, union of $(D)$-cells, which is normally radial (in the sense of section 1), and
has at most a singular point at the barycenter of each $(D)$-cell of top dimension $2m$, $m$ being
the complex dimension of $M$. This defines a cochain in the usual way, by assigning to each such cell
the Schwartz index of this vector field ({\it i.e.}  its Poincar\'e-Hopf index in the ambient space). This 
cochain is a cocycle that represents a cohomology class in   
$H^*(U, U\setminus V) \cong H^*(M, M\setminus V)$, and this class is by definiton the top Schwartz class
of $V$. Its image in $H_0(V)$ under Alexander duality gives the Euler-Poincar\'e characteristic
of $V$. The Schwartz classes of lower degrees are defined similarly, considering stratified $r$-frames,
$r = 1, \cdots, n = \rm{dim}\,V$, defined by radial extension on the $(2m-2r+2)$-skeleton of $(D)$, and
the corresponding  Schwartz index of such frames: just as the concept of Schwartz index  can
be extended to stratified vector fields in general (section 1), so too one can define the Schwartz index
of stratified frames in general, and use any such frame to define the corresponding Schwartz class
(see \cite {BLSS2, BSS3}).

The second  extension of the concept of tangent bundle is given by  the Nash bundle 
$\widetilde T \to  \widetilde V$ over the  Nash Transform $\widetilde V$,  which is 
 biholomorphic to $V_{reg}$ away from the divisor $\nu^{-1}(Sing(V))$, where 
$\nu: \widetilde V \to V$ is the projection.
Thus $\widetilde T$ can be regarded as a bundle that extends
 $T(V_{reg})$ to a bundle over $ \widetilde V$. The Chern classes of  $\widetilde T$ lie in
$H^*(\widetilde V)$, which is mapped into $H_*(\widetilde V)$ by the Alexander homomorphism
(see \cite {Br1, Br2});  the Mather classes of $V$, introduced in \cite {MP}, are by definition the
image of these classes under 
the morphism $\nu: H_*(\widetilde V) \to H_*( V)$.  MacPherson's
Chern classes for singular varieties \cite {MP},  lie in the homology of $V$ and 
can be thought of as being the Mather classes of $V$
weighted (in a sense that is made precise below) by the local Euler obstruction (\S 2 above). 
In fact, it is easy to show that the local Euler  obstruction 
satisfies that there exists unique integers
$\{n_i\}$ for which the equation 
\begin{equationth}\label{x.x}
\sum n_i \,Eu_{\overline V_\a}(x) = 1\, 
\end{equationth}
is satisfied 
 for all points $x$ in $V$, where the sum runs over all strata 
 $V_\a$ containing $x$ in their closure. Then  the  {\it MacPherson class} 
of degree $r$ is defined by:
$$
c_{r}(V)\, =\, c^M_{r}(\sum n_i \,{\overline V_\a}) \,=\, 
\sum n_i \,\iota_* c^M_{r}({\overline V_\a}) \,,$$
where $c^M_{r}({\overline V_\a})$ is the Mather class of degree $r$ of the analytic variety
${\overline V_\a}$.

 MacPherson's classes
satisfy important axioms and functoriality properties conjectured by Deligne and Grothendieck 
in the early 1970's.

Later, Brasselet and Schwartz  \cite {BS}   proved 
that the Alexander isomorphism
$H^*(M, M\setminus V) \cong H_*(V)$, carries the Schwartz classes into 
MacPherson's classes, so they are now called the {\it Schwartz-MacPherson classes} of $V$. As we briefly
explained before,  Schwartz classes are defined via the Schwartz indices of vector fields and frames;
the MacPherson classes are defined from the Chern classes of the Nash bundle (which determine the Mather 
classes) and the local Euler obstructions.

A third way of extending the concept of tangent bundle to singular varieties was introduced by 
Fulton and Johnson \cite {FJ}.
 The starting point is that if a variety $V \subset M $ is defined by a 
regular section $s$ of a 
holomorphic bundle $E$ over $M$, then one has   the virtual tangent 
 bundle $\tau V = [TM\vert_V - E|_V]$, introduced in \S 4 above.
 The Chern classes of the virtual tangent bundle 
 $\tau V$ (cap product the fundamental cycle $[V]$)  are the 
{\it Fulton-Johnson classes} of $V$. One has (\cite {LSS, SS3}) that the 0-degree Fulton-Johnson class
of such varieties equals the total virtual index of every continuous vector field with isolated 
singularities  on $V_{reg}$. Similar considerations hold for the higher degree Fulton-Johnson
classes, using frames and the corresponding virtual classes (see \cite {BLSS2, BSS3}).

Summarizing,  the various indices we presented in sections 1-4 are closely related to various
known characteristic classes of singular varieties that generalize the concept of Chern classes of 
complex manifolds. There is left the homological index of section 5:  this ought to be related with 
 a fourth  way for extending the concept of tangent bundle to singular varieties (with its corresponding
generalisation of Chern classes), introduced and studied by Suwa in \cite{Su6}. This is  by 
considering the
tangent sheaf $\Theta_V$, which is by definition the dual of $\Omega_V$,  the sheaf
of Kh\"aler differentials on $V$ introduced in \S 5. The latter is defined by the exact sequence:
$$ I_V/I^2_V \buildrel{f \mapsto df \otimes 1}\over{\longrightarrow}  \Omega_M \otimes \O_V \to  
\Omega_V  \to 0 \,,$$
and  $\Theta_V := \hbox{Hom}(\Omega_V, \C)$.
Both sheaves $\Omega_V$ and  $\Theta_V$ are coherent sheaves and 
 one can use them to define characteristic 
classes of $V$ that coincide with the usual 
 Chern classes when $V$ is non-singular. 
In particular, if $V$ is a local complete intersection in $M$, then one has a canonical locally free
resolution of $\Omega_V$ and the corresponding Chern classes essentially coincide with the 
Fulton-Johnson classes, though the corresponding classes for  $\Theta_V$ differ from these.
Recent work of J. Sch\"urmann  points out in this direction, at least if one considers the 
homological index of 1-forms.

We refer to \cite {Br2} for a rather complete presentation of characteristic classes of singular varieties, 
including the constructions of  Schwartz and MacPherson that we sketched above, and to \cite {BSS3} for
a discussion of indices of vector fields and their relation with 
characteristic classes of singular varieties, 
 much deeper than the one we presented here.

\bibliographystyle{plain}

$\,$

Jos\'e Seade

Instituto de Matem\'aticas, UNAM, 

Unidad Cuernavaca,

A. P. 273-3, Cuernavaca, Morelos, 

M\'exico

\end{document}